\documentstyle[amstex,amssymb, 11pt]{article}

\renewcommand{\epsilon}{\varepsilon}
\renewcommand{\phi}{\varphi}

\newtheorem{theorem}{Theorem}[section]

\newtheorem{lemma}[theorem]{Lemma}
\newtheorem{corollary}[theorem]{Corollary}

\newtheorem{example}[theorem]{Example}

\begin{document}
\thispagestyle{empty}

\begin{center}
{\Large\bf Herrero's Approximation Problem for Quasidiagonal Operators}
\end{center}

\begin{center}{\bf Nathanial P. Brown \footnote{ 
      Currently an MSRI Postdoctoral Fellow.}}\\
      UC-Berkeley\\
   Berkeley, California 94720 \\
   nbrown{\char'100}math.berkeley.edu 

\end{center}

\begin{abstract}
 Let $T$ be a quasidiagonal operator on a separable Hilbert space.  It
 is shown that there exists a sequence of operators $\{ T_n \}$ such
 that $dim( C^*(T_n)) < \infty$ and $\| T - T_n \| \to 0$ if and only
 if $C^*(T)$ is exact.
\end{abstract}

\vspace{4mm}

\begin{center}
{\it Dedicated to Yuriko Hoshiya and Isaac Seiji Brown.  }
\end{center}

\parskip2mm

\section{Introduction}

Let $H$ be a separable Hilbert space and $B(H)$ denote the bounded
operators on $H$.  $T \in B(H)$ is called block diagonal if there
exists an increasing sequence of finite rank projections $P_1 \leq P_2
\leq \ldots$ such that $P_n T - T P_n = 0$ for all $n$ and $P_n \to
1_H$ (s.o.t.).  $T \in B(H)$ is called {\em quasidiagonal} if it is
the norm limit of a sequence of block diagonal operators.

Herrero asked whether every quasidiagonal operator is the norm limit
of operators $T_n$ such that the $C^*$-algebra generated by each $T_n$
is finite dimensional.  This is equivalent to asking if quasidiagonal
operators are always norm limits of block diagonal operators each of
which is comprised of blocks of bounded dimension.  (See [Vo2,section
5] for a nice discussion of this problem, [DHS] for equivalent
formulations and [Vo2], [Br] for the basic theory of quasidiagonal
(sets of) operators.)  A negative answer was first obtained by
S. Szarek (cf.\ [Sz]).  Using probabalistic methods he was able to
show the existence of a quasidiagonal operator for which there was an
operator {\em theoretic} obstruction to such approximations.  In [Vo1]
Voiculescu observed an operator {\em algebraic} obstruction to the
existence of such approximations which we now recall.

Assume $T_n \to T$ (in norm) and $dim( C^*(T_n)) < \infty$ for all
$n$.  Let $E_n : B(H) \to C^*(T_n)$ be conditional expectations (which
exist by finite dimensionality).  It follows that $E_n (x) \to x$ (in
norm) for every $x \in C^*(T)$ and hence the inclusion $C^*(T)
\hookrightarrow B(H)$ is a nuclear map.  Thanks to the work of
E. Kirchberg, this latter condition is equivalent to saying that
$C^*(T)$ is an exact $C^*$-algebra (cf.\ [Wa]).  Hence exactness
provides a natural operator algebraic obstruction to Herrero's
approximation question.  In this note we observe that this is the only
obstruction.

\begin{theorem}
If $T$ is a quasidiagonal operator on a separable Hilbert space, 
 then $T$ is the norm limit of operators $\{ T_n \}$ such that $dim(
 C^*(T_n)) < \infty$ for every $n$ if and only if $C^*(T)$ is exact.
\end{theorem}

This gives an affirmative answer to [DHS, Problem 4.3] and is an
immediate consequence of the following answer to a question of
Dadarlat (cf.\ [Da]).

\begin{theorem}
Let $A$ be an exact $C^*$-algebra and $\psi : A \to B(H)$ be a
*-homomorphism with $H$ separable.  Then $\psi(A)$ is a quasidiagonal
set of operators if and only if for each $\epsilon > 0$ and finite set
${\cal F} \subset A$ there exists a finite dimensional
$C^*$-subalgebra $C \subset B(H)$ with $\psi({\cal F})
\subset^{\epsilon} C$ (i.e.\ for each $a \in {\cal F}$ there exists $c
\in C$ with $\| \psi(a) - c \| < \epsilon$).
\end{theorem}

Dadarlat proved Theorem 1.2 under the additional assumption that
$\psi(A)$ contains no nonzero compact operators (cf.\ [Da, Thm.\ 6]).
We will deduce the general case from Dadarlat's result and a technical
generalization of Voiculescu's noncommutative Weyl-von Neumann
Theorem.

\section{Proofs of Main Results}

If $A$ is a separable, unital $C^*$-algebra and $\phi : A \to B(H)$
(with $H$ separable and infinite dimensional) is a unital completely
positive map then we say that $\phi$ is a {\em faithful representation
modulo the compacts} if $\pi \circ \phi : A \to Q(H)$ is a
*-monomorphism, where $\pi$ is the quotient map onto the Calkin
algebra.  In this situation we define constants $\eta_{\phi} (a)$ by
$$\eta_{\phi} (a) = 2\max( \| \phi(a^* a) - \phi(a^*) \phi(a)
\|^{1/2}, \| \phi(a a^*) - \phi(a) \phi(a^*) \|^{1/2} )$$ for every $a
\in A$.  The following lemma is the generalization of Voiculescu's
Theorem that we will need.  The main idea in the proof is essentially
due to Salinas (see the proofs of [Sa, Thm.\ 2.9] and [DHS, Thm.\
4.2]).  This result also appears in [Br], but we include the proof for
the reader's convenience.

\begin{lemma}
Let $A$ be a separable, unital $C^*$-algebra and $\phi : A \to B(H)$
be a faithful representation modulo the compacts.  If $\sigma : A \to
B(K)$ is any faithful, unital, representation such that $\sigma(A)$
contains no nonzero compact operators then there exist unitaries $U_n
: H \to K$ such that $$\limsup\limits_{n \to \infty} \| \sigma(a) -
U_n \phi(a) U_n^* \| \leq \eta_{\phi} (a)$$ for every $a \in A$.
\end{lemma}

{\noindent\it Proof.} Note that by the usual version of Voiculescu's
Theorem it suffices to show that {\em there exists} a representation
$\sigma$ satisfying the conclusion of the theorem since all such
representations are approximately unitarily equivalent.

Let $\rho : A \to B(L)$ be the Stinespring dilation of $\phi$; i.e.\
$\rho$ is a unital representation of $A$ and there exists an isometry
$V : H \to L$ such that $\phi(a) = V^* \rho(a) V$, for all $a \in A$.
Let $P = VV^* \in B(L)$ and $P^{\perp} = 1_L - P$.  A routine
calculation shows that for every $a \in A$, $$\| P^{\perp}\rho(a) P \|
\leq \| \phi(a^* a) - \phi(a^*) \phi(a) \|^{1/2}.$$

Now write $L = PL \oplus P^{\perp}L$ and decompose the
representation $\rho$ accordingly.  That is, consider the matrix
decomposition $$\rho (a) =
\left(\begin{array}{cc}
\rho (a)_{11} & \rho (a)_{12} \\
\rho (a)_{21} & \rho (a)_{22}
\end{array}\right), $$  
where $\rho (a)_{21} = P^{\perp}\rho(a) P$ and $\rho (a)_{12} = \rho
(a^*)_{21}^*$.    

Consider the Hilbert space $P^{\perp}L \oplus
PL$ and the representation $\rho^{\prime} : A \to B(P^{\perp}L \oplus
PL)$ given in matrix form as
$$\rho^{\prime} (a) =
\left(\begin{array}{cc}
\rho (a)_{22} & \rho (a)_{21} \\
\rho (a)_{12} & \rho (a)_{11}
\end{array}\right).$$ 
Now using the obvious identification of the Hilbert spaces $$PL \oplus
\big(\bigoplus_{\Bbb N} P^{\perp}L \oplus PL \big) \ \ {\rm and} \ \
\bigoplus_{\Bbb N} L = \bigoplus_{\Bbb N}( PL \oplus P^{\perp}L)$$ a
standard calculation shows that $$\| \rho(a)_{11} \oplus
\rho^{\prime\infty}(a) - \rho^{\infty}(a)\| \leq \eta_{\phi}(a)$$ for
all $a \in A$, where $\rho^{\prime\infty} = \oplus_{\Bbb N}
\rho^{\prime}$ and $\rho^{\infty} = \oplus_{\Bbb N}
\rho$.  Note also that $\rho(a)_{11} = V \phi(a) V^*$.

Now, let $C$ be the linear space $\phi (A) + {\cal K}(H)$.  Note that
$C$ is actually a separable, unital $C^*$-subalgebra of $B(H)$ with
$\pi (C) = A$ where $\pi: B(H) \to Q(H)$ is the quotient map onto the
Calkin algebra.  By Voiculescu's Theorem we have that $\iota \oplus
\rho^{\prime\infty} \circ \pi$ is approximately unitarily equivalent
modulo the compacts to $\iota$, where $\iota : C \hookrightarrow B(H)$
is the inclusion.  Let $W_n : H \to H \oplus (\oplus_{\Bbb N}
(P^{\perp}L \oplus PL))$ be unitaries such that $$\| \phi (a) \oplus
\rho^{\prime\infty} (a) - W_n \phi (a) W_n^* \| \to 0$$ for all $a \in
A$.

Let $$\tilde{V} : H \oplus (\bigoplus_{\Bbb N} (P^{\perp}L
\oplus PL )) \to \bigoplus_{\Bbb N} L$$ be the unitary $V \oplus 1$
(again using the obvious identification of $PL \oplus
\big(\oplus_{\Bbb N} (P^{\perp}L \oplus PL) \big)$ and $\oplus_{\Bbb
N} L$).  Note that $\tilde{V}(\phi(a) \oplus
\rho^{\prime\infty}(a))\tilde{V}^* = V\phi(a)V^* \oplus
\rho^{\prime\infty}(a) = \rho(a)_{11} \oplus
\rho^{\prime\infty}(a)$. We complete the proof by defining $$K =
\oplus_{\Bbb N} L, \ \ \ \sigma = \rho^{\infty} = \oplus_{\Bbb N}
\rho, \ \ \ U_n = \tilde{V}W_n : H \to \oplus_{\Bbb N} L = K. \ \ \
\Box$$

{\noindent\bf Proof of Theorem 1.2.}  We may assume that $A$ is
separable unital and $\pi$ is a unital *-homomorphism. Since quotients
of exact $C^*$-algebras are again exact it suffices to prove that if
$A \subset B(H)$ is an exact $C^*$-algebra and a quasidiagonal set of
operators then for each finite set ${\cal F} \subset A$ and $\epsilon
> 0$ there exists a finite dimensional subalgebra $C \subset B(H)$
such that ${\cal F} \subset^{\epsilon} C$.  In fact, it suffices to
show that ${\cal F} \subset^{\epsilon} C$ where $C$ is an
approximately finite dimensional (AF) $C^*$-subalgebra of $B(H)$.

Let $P_1 \leq P_2 \leq P_3 \ldots$ be an increasing sequence of finite
rank projections which converge to the identity operator in the strong
operator topology and such that $\| aP_n - P_n a \| \to 0$ for all $a
\in A$.  Let $B \subset Q(H)$ be the image of $A$ in the Calkin
algebra.  Then we have the natural exact sequence $0 \to {\cal K}(H)
\to A + {\cal K}(H) \to B \to 0$.  Since ${\cal K}(H)$ and $B$ are
exact and the quotient map $A + {\cal K}(H) \to B$ is locally liftable
(since $A$ is exact), it follows from [Ki2, Prop.\ 7.1] that $A +
{\cal K}(H)$ is also exact. Since the compact operators are actually
nuclear it follows from [EH], that this short exact sequence is
semisplit; i.e.\ there exists a unital completely positive splitting
$\phi : B \to A + {\cal K}(H)$.  Define $\phi_n (a) = (1 - P_n)
\phi(a) (1 - P_n)$ and note that each $\phi_n$ is also a completely
positive splitting.  Moreover (as is well known) the maps $\phi_n$ are
asymptotically multiplicative and hence $B$ is a quasidiagonal (QD)
$C^*$-algebra. Let $H_n = (1 - P_n)H$. Note that each $H_n$ has finite
codimension in $H$.

Let $\rho : B \to
B(K)$ be a faithful representation such that $\rho(B)$ contains no
nonzero compact operators.  By Dadarlat's result ([Da, Thm.\ 6]), we
can find a sequence of finite dimensional $C^*$-subalgebras $D_n
\subset B(K)$ such that each element of $\rho(B)$ is the norm limit of
a sequence taken from the algebras $D_n$.  Since the maps $\phi_n$ are
asymptotically multiplicative faithful representations modulo the
compacts (when regarded as taking values in $B(H_n)$) and each $H_n$
has finite codimension in $H$ it follows from Lemma 2.1 that we can
find a sequence of finite dimensional $C^*$-subalgebras $E_n \subset
Q(H)$ such that each element of $B$ is the norm limit of a sequence
taken from the algebras $E_n$.  Since extensions of AF algebras are
again AF the proof of the theorem is completed by defining $C_n$ to be
the pullbacks of the $E_n$.  $\Box$

In [Sa] N. Salinas defined $Ext_{qd} (B)$ for a QD $C^*$-algebra as
follows.  A unital *-monomorphism $\rho : B \to Q(H)$ defines an
element of $Ext_{qd} (B)$ if the pullback of $\rho(B)$ in $B(H)$ is a
quasidiagonal set of operators.  It is easy to see that this notion in
invariant under strong (unitary) equivalence of extensions (i.e.\
gives a well defined subset of $Ext(B)$).  When $B$ is nuclear,
$Ext_{qd} (B)$ is a subgroup of $Ext(B)$ (cf.\ [Sa, Cor.\ 2.10]).
Kirchberg has given examples of exact QD $C^*$-algebras $B$ and
*-monomorphisms $\rho : B \to Q(H)$ which define noninvertible
elements of $Ext_{qd} (B)$ (cf.\ [Ki1]).  The proof of Theorem 1.2
above easily yields a new characterization (different than Arveson's
characterization as those *-monomorphisms which admit a completely
positive lifting) of the invertible elements of $Ext_{qd} (B)$
when $B$ is exact.

\begin{corollary}
Let $B$ be a separable unital exact QD $C^*$-algebra and $\rho : B \to
Q(H)$ be a *-monomorphism such that $[\rho] \in Ext_{qd} (B)$.  Then
$[\rho]$ is invertible if and only if for each finite subset ${\cal F}
\subset B$ and $\epsilon > 0$ there exists a finite dimensional
$C^*$-subalgebra $C \subset Q(H)$ such that $\rho( {\cal F} )
\subset^{\epsilon} C$.
\end{corollary}
 
{\noindent\bf Proof.}  If $[\rho]$ is invertible then there exists a
completely positive splitting and hence the proof of Theorem 1.2 shows
that $\rho(B)$ admits such approximants.  On the other hand, if
$\rho(B)$ admits such approximants then the argument in the
introduction shows that $\rho$ is a nuclear map and hence, by the
Choi-Effros Lifting Theorem, is liftable (hence invertible).  $\Box$

We wish to thank M. Dadarlat for pointing out the following  example.  

\begin{example}
{\em Let $H$ be a Hilbert space with orthonormal basis $\{ e_i \}_{i
\in {\Bbb Z}}$.  Let $U \in B(H)$ be the bilateral shift unitary
(i.e.\ $U(e_i) = e_{i + 1}$) and $T \in B(H)$ be any operator with the
property that each $e_i$ is an eigenvector for $T$.  If we define $S =
TU$ then $C^*(S)$ is exact.  This is because it is contained in a
nuclear $C^*$-algebra.  More precisely, the $C^*$-algebra generated by
$U$ and the set of {\em all} such $T$ (as above) is naturally
isomorphic to the crossed product $C^*$-algebra $l^{\infty}({\Bbb Z})
\rtimes_{\gamma} {\Bbb Z}$ where the action $\gamma$ on
$l^{\infty}({\Bbb Z})$ arises from the bilateral shift $i \mapsto i +
1$ on ${\Bbb Z}$.  Since crossed products of nuclear $C^*$-algebras by
amenable groups are again nuclear, this shows that $C^*(S) \subset
l^{\infty}({\Bbb Z}) \rtimes_{\gamma} {\Bbb Z}$ is an exact
$C^*$-algebra.  Evidently this example extends to one sided (i.e.\
unilateral) shifts as well since these correspond to crossed products
by endomorphisms (which still preserves nuclearity).}
\end{example}

The example above is nothing but the class of weighted shift (both
bilateral and unilateral) operators.  R. Smucker first characterized
the quasidiagonality of weighted shifts in terms of the associated
weight sequence in [Sm].  Thus our main result gives a new proof of a
result of D. Herrero who showed that all quasidiagonal weighted shifts
admit Herrero approximations (see [He]).  Note, however, that the
present approach has the added benefit of treating $n$-tuples as well.
That is, our results show that if $S_1, \ldots, S_k$ are weighted
shifts (with respect to the same basis) and $\{ S_1, \ldots, S_k \}$
is a quasidiagonal {\em set} of operators then for each $\epsilon > 0$
there exist operators $\tilde{S}_1, \ldots, \tilde{S}_k$ such that $\|
S_i - \tilde{S}_i \| < \epsilon$ and $C^*(\tilde{S}_1, \ldots,
\tilde{S}_k)$ is finite dimensional (since $C^*( S_1, \ldots, S_k)$
will always be exact).

\end{document}